\documentclass[10pt]{article}
\usepackage{amsmath,amssymb,amsthm,amscd}
\numberwithin{equation}{section}

\def\p{\partial}

\newtheorem{prop}{Proposition}[section]
\newtheorem{theo}[prop]{Theorem}
\newtheorem{lem}[prop]{Lemma}
\newtheorem{cor}[prop]{Corollary}
\newtheorem{rem}[prop]{Remark}

\newtheorem{defi}[prop]{Definition}
\newtheorem{conj}[prop]{Conjecture}

\def\begeq{\begin{equation}}
\def\endeq{\end{equation}}

\def\and{\quad{\rm and}\quad}

\let\lra=\longrightarrow

\def\mapright\#1{\,\smash{\mathop{\lra}\limits^{\#1}}\,}

\begin {document}

\bibliographystyle{plain}
\title{A new parabolic flow in K\"ahler manifolds}
\author{X. X. Chen }
\date{Revised on Sep. 26, 2000}
\maketitle


\subsection{Introduction and main theorem}
This is a continuation of my earlier work \cite{chen993}, where
we study the lower bound of the Mabuchi energy on
 K\"ahler manifolds when the first Chern class is negative. The
 problem of finding a lower bound for the Mabuchi energy is very
important in K\"ahler geometry since the existence of a lower
bound is the pre-condition for the existence of constant scalar
curvature metric in a K\"ahler class (cf. \cite{Bando87} and
\cite{chen991} ).  The Mabuchi energy is first defined by Mabuchi
through its derivative. The first explicit formula was given in
\cite{tian98}.  In \cite{chen993}, we re-group the right hand side
of this formula of Tian in a slightly different way. According to
this decomposition formula (cf. \cite{chen993}), the problem is
reduced to the solution of the existence problem for critical
metrics of a new functional $J_{\chi}$ (cf. \cite{chen993}
\cite{Dona99}). For convenience, we include the definition of the
functional $J$ below. The existence problem is completely solved
in any K\"ahler surface. However, it is still widely open  in
higher dimensional K\"ahler manifolds.  In this paper, we try to
understand the general existence problem via the flow method. Let
$(V^n, \omega_0)$ be an $n$ dimensional K\"ahler manifold and
$\omega_0$ be a fixed K\"ahler form in $V.\;$ Consider the space
of K\"ahler potentials
\[{\cal H} = \{\varphi\in C^{\infty}(V) \mid\; \omega_{\varphi} = \omega_0 + i \partial \overline{\partial} \varphi > 0 \;{\rm on}\;V \}.
\]
For any positive (1,1) form $\chi$,  consider the parabolic
equation
\begin{equation}
{{\partial \varphi}\over {\partial t}} = {{n\;[\chi] \cdot
[\omega_0]^{n-1}} \over {[\omega_0]^n}} - {{\chi \wedge
{\omega_{\varphi}}^{n-1}} \over \omega_{\varphi}^n} .
\end{equation}
 This is the gradient flow for the
functional $J_\chi$ which will be defined in Section 0.2. Similar
to the case of the Calabi flow \cite{chen992}, the main result of
this paper is

\begin{theo}  The following statements are true:
\begin{enumerate}
\item
This gradient flow of $J_{\chi}$ always exists for all the time
for
 any smooth initial data.
Moreover, the length of any smooth curve in $\cal H$  and  the
distance between any two metrics decreases under this
flow\footnote{The definition of the length of a smooth curve in
$\cal H$ is given in equation (\ref{eq:length}). On the other
hand, the distance function in $\cal H$ is given by Definition
1.3. }.
\item If the bisectional curvature of $\chi$ is semi-positive,
then the gradient flow exists for all time and converges to a
smooth critical metric.
\end{enumerate}

\end{theo}

In this paper, we will call the gradient flow of the functional
$J_\chi$ a $J-$flow.

\subsection{Setup of notations}

Let $g_0 = \displaystyle \sum_{\alpha,\overline{\beta}=1}^n
{g_0}_{\alpha,\overline{\beta}}  d\;z_{\alpha} d\;{\bar
z}_{\beta}$ be the K\"ahler metric corresponding to the K\"ahler
form $\omega_0 $ and $g = \displaystyle
\sum_{\alpha,\overline{\beta}=1}^n {g}_{\alpha,\overline{\beta}}
d\;z_{\alpha} d\;{\bar z}_{\beta}$   the K\"ahler metric
corresponds to the K\"ahler form $\omega_{\varphi}$ for some
$\varphi \in \cal H.\;$      For any fixed positive (1,1) form
$\chi, $ one introduces a new functional $J_\chi$ with respect to
this form $\chi.\;$

\begin{defi}Suppose $\chi$ is a positive closed (1,1) form in $V.\;$ For any
$\varphi(t) \in \cal H,$  the functional $J_{\chi}$ is defined
through its derivative:
\[
{{d\,J_{\chi}}\over{d\,t}} = \int_V \;{{\partial
\varphi}\over{\partial t}} \chi\wedge
{{{\omega_{\varphi}}^{n-1}}\over {(n-1)!}}.
\]
It is straightforward to show that this is well defined.
\end{defi}
\begin{rem} This version of definition is given by \cite{Dona99}. In
\cite{chen993}, we give an independent definition for the
functional $J_{\chi}$ where we assume $\chi$ is a Ricci form (
not necessary positive).
\end{rem}

 In \cite{Dona99},
Donaldson pointed out the significance of this functional in its
own right: a) $J_{\chi}$ is a convex functional in the space of
K\"ahler potentials (See proposition 3.1 below); b)  $J_{\chi}$
is a moment map from the space of K\"ahler potentials to  the
 dual space of the Lie algebra of some sympletic automorphism group.
  \\

For any $\varphi\in \cal H,$ there exists a unique K\"ahler form
$\omega_{\varphi}.\; $ Conversely, for  any K\"ahler form $\omega'
$ in the cohomology class of $\omega_0,$ there exists a unique
K\"ahler potential $\varphi\in \cal H $ up to additiion of some
constants such that
\[
\omega' = \omega_{\varphi}.
\]

In order to remove the ambiguity of the correspondence between
K\"ahler forms and K\"ahler potential, we introduce a well-known
functional here.
\begin{defi} For any smooth curve $\varphi(t) \in \cal H,$ we
define the functional $I$ such as
\[
   {{d\, I(\varphi(t))}\over {d\,t}} = \displaystyle \int_V\; {{\p
   \varphi}\over {\p t}}\; {{\omega_{\varphi}^n}\over {n!}}.
\]
Clearly, $I$ is well defined.
\end{defi}

Obviously, $I$ defines a function from $\cal H$ to the real line.
Consider the $0-$level surface ${\cal H}_0 $ of the functional
$I$ in ${\cal H}.\;$ This level surface  can be defined as
  \[
  {\cal H}_0 =\{\varphi \in {\cal H}\;\mid\; I( \varphi)
  =0\}.
  \]

  Set the variational space of $J_{\chi}$ as ${\cal H}_0.\;$ In local coordinates, suppose that $\chi$ can be expressed as
  the following:
\[
\chi =  \displaystyle \sum_{\alpha, \beta=1}^n\; \chi_{\alpha
\overline{\beta}}\; d\;z_{\alpha} d\;{\bar z}_{\beta}.
\]
    Then the Euler-Lagrange  equation for $J_{\chi}$
   is:
\begin{equation}
  \displaystyle \sum_{\alpha, \beta=1}^n\; g^{\alpha \overline{\beta}}\; \chi_{\alpha \overline{\beta}}
  =    tr_{g} \chi = c
\label{eq:Euler}
  \end{equation}
where
 \begin{equation}
   g_{\alpha \bar \beta} = {g_0}_{\alpha \bar \beta} + {{\p^2 \varphi}\over {\p z^{\alpha} \bar \p z^{\beta}}},
\label{eq:evolvedmetric}
 \end{equation}
Here $g = \displaystyle \sum_{\alpha,\beta=1}^n\; g_{\alpha\bar
\beta} d\,z^{\alpha} \;d\,z^{\bar \beta} $ is the critical metric
of $J_\chi $ in ${\cal H}_0$ and $c$ is a constant which depends
only on the K\"{a}hler class of $[\chi]$ and $[\omega_0], i.e.,$
\begin{equation}
c = {{\displaystyle \int_V \; \chi \wedge {1\over {(n-1)!}}\omega^{(n-1)}}\over{ \displaystyle \int_V \; {1 \over {n!}} \omega^n}}
 =  {{n\;[\chi] \cdot [\omega_0]^{n-1}} \over {[\omega_0]^n}}.
\label{eq:const}
\end{equation}

From equation (\ref{eq:Euler}), it is easy to see that a
necessary condition for a solution to exist is (also see
\cite{Dona99})
\[
  c \cdot \omega_g - \chi > 0,
\]
where $\omega_g$ is the K\"{a}hler form associated to  $g.\;$
In other words, there exists at least one K\"{a}hler form $\omega$ in the K\"{a}hler
class of $[\omega_0]$ such that the following holds:
\begin{equation}
  c \cdot \omega - \chi > 0.
\label{eq:necessary1}
\end{equation}

\begin{conj} (Donaldson \cite{Dona99}) If the aforementioned
necessary condition is satisfied, then there exists a critical
 point for $J_{\chi}$ in that
K\"ahler class.
\end{conj}

\noindent {\bf Historic remarks}: Using the heat flow method to
study the nonlinear PDE is a well known method in differential
geometry. In recent years, it has been the source of active
researches since the famous work of J. Eells  and J. Sampson
\cite{ES64}. Interested readers are refereed to important works by
R. Hamilton  \cite{Hamilton93}, G. Huisken
\cite{Huisken90}\cite{Huisken98} and a survey paper by H. D. Cao
 and B. Chow \cite{Cao99} and the
references therein.\\

{\bf Acknowledgement} The author would like to thank Professor
Donaldson for his encouragement in studying the critical points
for the functional $J.\;$  He also wants to thank Professor
 G. Huisken  for many interesting discussion about heat flow. He also wants  to thank
R. Schoen, L. Simon, E. Calabi  and G. Tian for their interests
in this work. Thanks also to  Guofang Wang for carefully reading
an earlier version of this paper. Finally, thanks also to the
referees for many interesting suggestions.

\section{Summary of recent developments in the Riemannian metric
in space of K\"ahler potentials}
 Mabuchi (\cite{Ma87})  in 1987 defined a Riemannian metric
on the space of K\"ahler metrics,
under which it  becomes (formally) a non-positive curved infinite dimensional
 symmetric space. Apparently unaware of Mabuchi's work,
Semmes \cite{Semmes92}  and Donaldson \cite{Dona96}  re-discovered
this same metric  from different angles. For any  vector $\psi$
in the tangential space $ T_{\varphi} \cal {H}, $ we define the
length of this vector as
\[
\|\psi\|^2_{\varphi} =\int_{V}\psi^2\;d\;\omega_{\varphi}^n.
\]
For any smooth curve $\varphi(t): [0,1] \rightarrow \cal H,$ the
length of this curve is given by \begin{equation} \displaystyle
\int_0^1 \left( \displaystyle \int_V\; \left({{\p
\varphi(t)}\over {\p t}}\right)^2 \omega_{\varphi}^n
\right)^{1\over 2 }\;d\,t. \label{eq:length}
\end{equation}
 The geodesic equation is
\begin{equation}
  {{\p^2 \varphi(t)}\over {\p
t^2}} - {1\over 2} |\nabla {{\p \varphi(t)}\over {\p t}}|^2_{
\varphi(t)} = 0, \label{geodesic}
\end{equation}
where the derivative and norm in the 2nd term of the left hand side
are taken with respect to the metric $\omega_{\varphi(t)}.\;$\\

This geodesic equation shows us how to define a connection on the
tangent bundle of ${\cal H}$.  If $\phi(t)$ is any path in ${\cal
H}$ and $\psi(t)$ is a field of tangent vectors along the path
(that is, a function on $V \times [0,1]$), we define the
covariant derivative along the path to be
$$ D_{t}\psi = \frac{\partial\psi}{\partial t} - {1\over 2}  (\nabla \psi, \nabla {{\p \phi(t)}\over {\p
t}})_{\phi}.  $$
 The main theorem formally
proved  in \cite{Ma87} (and later reproved in \cite{Semmes92} and
 \cite{Dona96}) is:\\

\noindent {\bf Theorem A} {\it The Riemannian manifold $\cal {H} $ is an infinite dimensional symmetric space; it admits a Levi-Civita connection whose
curvature is covariant constant. At a point $\phi\in{\cal {H}}$ the curvature  is given by
\[     R_{\phi}(\delta_{1}\phi, \delta_{2}\phi) \delta_{3}\phi=
- {1\over 4}  \{ \{ \delta_{1}\phi, \delta_{2}\phi\}_{\phi},
\delta_{3}\phi\}_{\phi},\]
where $\{\ ,\ \}_{\phi}$ is the Poisson bracket on $C^{\infty}(V)$ of the
symplectic form $\omega_{\phi}$; and $\delta_1 \phi, \delta_2 \phi \in T_{\phi} {\cal H}.\;$ Then the sectional curvature is non-positive, given by}
\[    K_{\phi}(\delta_{1}\phi, \delta_{2}\phi) = - {1\over 4}  \Vert
 \{ \delta_{1} \phi , \delta_{2}\phi\}_{\phi}\Vert_{\phi}^2. \]
We will skip the proof here. Interested readers are referred to
paper of Mabuchi \cite{Ma87} or \cite{Semmes92} and
 \cite{Dona96} for the proof.\\

 In \cite{Dona96}, Donaldson  outlined a connection between
this Riemannian metric in the infinite dimensional space $\cal H$
and the traditional K\"ahler geometry through a series important
conjectures and theorems. In 1997, following his program, the
author proved some of
his conjectures:\\

\noindent {\bf Theorem B} \cite{chen991}{\it The following statements
are true:
\begin{enumerate}
\item The space of K\"ahler potentials ${\cal H}$ is convex by $C^{1,1}$ geodesic (cf. Definition 1.1 below).
More specifically, if $\varphi_0,\varphi_1 \in \cal H$ and
$\varphi(t) \;(0\leq t \leq 1)$ is a geodesic connecting these
two points in $\cal H,$ then the mixed covariant derivative of
$\varphi(t)$ is uniformly bounded from above.
\item  ${\cal H}$ is a metric space. In other words, the infimum of the
lengths of all possible curves between any two points in $\cal H$ is strictly
positive.
\end{enumerate}
}

\begin{defi} For any $\epsilon > 0,$
a smooth path $\varphi(t)$ in $\cal H$ is called
$\epsilon$-approximate geodesic if the following holds:
\begin{equation}
\left({{\p^2 \varphi}\over {\p t^2}}  - {1\over 2} |\nabla {{\p
\varphi}\over {\p t}} |_{\varphi(t)}^2\right) \omega_{\varphi}^n
= \epsilon \cdot \omega_0^n. \label{eq:epsilongeodesic}
\end{equation}
\end{defi}

\begin{rem} Note that for $\epsilon \rightarrow 0,$ the
$\epsilon-$approximated geodesic converges to the unique
$C^{1,1}$ geodesic between the two end points. Theorem B implies
that the second mix derivatives of $\varphi(t)$ stay uniformly
bounded as $\epsilon \rightarrow 0.$
\end{rem}

 Following Theorem B, we can define a distance function in $\cal
H.\;$
\begin{defi} For any two points $\varphi_1,\; \varphi_2 \in \cal
H,\;$ the distance $d(\varphi_1,\varphi_2)$ is defined to be the
length of the unique $C^{1,1}$ geodesic which connects these two
points.
\end{defi}

In \cite{chen992}, E. Calabi and the author proved the following:\\

\noindent {\bf Theorem C}\cite{chen992}{\it The following statements are true:
\begin{enumerate}
\item  $\cal H$ is a non-positive curved space in the sense of
Alenxandrov.
\item  The length of any curve in $\cal H$ is decreased under
the Calabi flow unless it is represented by a holomorphic
transformation. The distance in $\cal H$ is also decreasing if
the Calabi flow exists for all the time for any initial smooth
data.
\end{enumerate}}

\section{General preparation}
 For convenience, from now on,  we drop the dependence of $\chi$ in
$J_\chi.\;$ Then
\begin{prop} $J$ is a strictly convex functional on any $C^{1,1}$
geodesic. In particular, $J$ has at most one critical point in
${\cal H}_0.\;$
\end{prop}
\begin{rem} This proposition was pointed out to the author by Donladson in
spring 1999. The proof here is somewhat different from his
original idea.
\end{rem}
\begin{proof} Suppose $\varphi(t)$  is a $C^{1.1}$ geodesic. In other
words, $\varphi(t)$ is a weak limit of the following continuous
equation as $\epsilon \rightarrow 0 $ (with uniform bounds on the
second mixed derivatives of K\"ahler potentials):
\[  \left({{\partial^2 \varphi}\over{\partial t^2}}  - {1\over 2} \;\mid \nabla {{\partial \varphi}\over{\partial t}} \mid_{\varphi}^2\right)
{{\omega_{\varphi(t)}^n}\over {n!}}  = \epsilon \cdot
{{\omega_0^n}\over {n!}} .
\]
Denote $g(t)$ is the corresponding K\"ahler metric corresponds to
the K\"ahler form $\omega_{\varphi(t)}.\;$  Again, we drop the
dependence of $t$ from $g(t)$ for convenience from now on. Recall
the  definition of $J$, we have
\[
   {{d\, J}\over {d\, t}} = \displaystyle\;\int_V\; {{\partial \varphi}\over{\partial t}}
   \; (g^{\alpha \overline{\beta}} \; \chi_{\alpha \overline{\beta}}) {{\omega_{\varphi}^n}\over {n!}} .\;
\]
Then (denote $\sigma = g^{\alpha \overline{\beta}} \; \chi_{\alpha \overline{\beta}}$ in the following calculation):
\[\begin{array}{lcl}
 & & {{d\,^2 J}\over {d\, t^2}} \\
 &  = & \displaystyle\;\int_V\;\left({{\partial^2 \varphi}\over{\partial t^2}} \sigma -  {{\partial \varphi}\over{\partial t}} g^{\alpha \overline{\beta}} ({{\partial \varphi}\over{\partial t}})_{,\overline{\beta} r} g^{r \overline{\delta}} \chi_{\alpha\overline{\delta}} +{{\partial \varphi}\over{\partial t}} \;\sigma \;\triangle_g\; {{\partial \varphi}\over{\partial t}}
\right) \; {{\omega_{\varphi}^n}\over {n!}} \\
& = &  \displaystyle\;\int_V\;\left( {{\partial^2
\varphi}\over{\partial t^2}} \sigma - {{\partial
\varphi}\over{\partial t}} g^{\alpha \overline{\beta}}
({{\partial \varphi}\over{\partial t}})_{,\overline{\beta} r}
g^{r \overline{\delta}} \chi_{\alpha\overline{\delta}}\right.
\\ & & \qquad \qquad \qquad-\left. ({{\partial \varphi}\over{\partial t}})_{,r} \sigma g^{r \overline{\delta}} ({{\partial \varphi}\over{\partial t}})_{,\overline{\delta}}
- {{\partial \varphi}\over{\partial t}} g^{\alpha \overline{\beta}} \chi_{\alpha \overline{\beta},\overline{\delta}}
g^{r \overline{\delta}} ({{\partial \varphi}\over{\partial t}})_{,r}  \right) \; {{\omega_{\varphi}^n}\over {n!}} \\
& = &  \displaystyle\;\int_V\;\left( ({{\partial^2 \varphi}\over{\partial t^2}}
- {1\over 2} |\nabla {{\partial \varphi}\over{\partial t}}|^2_{g} ) \sigma
- {{\partial \varphi}\over{\partial t}} g^{\alpha \overline{\beta}} ({{\partial \varphi}\over{\partial t}})_{,\overline{\beta} r}
 g^{r \overline{\delta}} \chi_{\alpha\overline{\delta}} \right. \\
 & & \qquad \qquad \qquad \left.- {{\partial \varphi}\over{\partial t}} \left(g^{\alpha \overline{\beta}} \chi_{\alpha \overline{\delta}} g^{r \overline{\delta}}
\right)_{,\overline{\beta}} ({{\partial \varphi}\over{\partial t}})_{,r} \right) \; {{\omega_{\varphi}^n}\over {n!}} \\
& = &  \displaystyle\;\int_V\;\left( ({{\partial^2 \varphi}\over{\partial t^2}}
- {1\over 2} |\nabla {{\partial \varphi}\over{\partial t}}|^2_{g} )\; ( g^{\alpha \overline{\beta}} \; \chi_{\alpha \overline{\beta}} )
+ ({{\partial \varphi}\over{\partial t}})_{,\overline{\beta}} \left(g^{\alpha \overline{\beta}} \chi_{\alpha \overline{\delta}} g^{r \overline{\delta}}
\right) ({{\partial \varphi}\over{\partial t}})_{,r} \right) \; {{\omega_{\varphi}^n}\over {n!}} \\
& \geq & \displaystyle\;\int_V\;({{\partial \varphi}\over{\partial
t}})_{,\overline{\beta}} \left(g^{\alpha \overline{\beta}}
\chi_{\alpha \overline{\delta}} g^{r \overline{\delta}} \right)
({{\partial \varphi}\over{\partial t}})_{,r} \;
{{\omega_{\varphi}^n}\over {n!}} \geq 0.
\end{array}
\]
The last equality holds along any $C^{1,1}$ geodesic.
\end{proof}

\begin{defi} According to \cite{Dona99}, the functional $J$ can be regarded as a moment map. Therefore,
we can define its "norm" as  a new energy  functional $E.\;$ In
other words, we have,
 \[E = \displaystyle
\int_V \; (tr_g (\chi))^2 \; {{\omega_{\varphi}^n}\over {n!}} .\]
\end{defi}
\begin{prop} $E$ has similar critical points as $J.\;$ Moreover, $E$ is decreasing along the
gradient flow of $J.\;$
\end{prop}
\begin{proof}
 Set $\sigma = tr_g \chi.\;$ All of the derivatives, norm and
integration  are taken with respect to $g$ in the following
calculation:
\[
\begin{array}{lcl}
  & & \delta_v \; E(g)\\
   & = & \displaystyle \int_V  \left( 2 \sigma (- \displaystyle \sum_{\alpha,\beta,r,\delta=1}^n g^{\alpha \overline{\beta}} v_{,\overline{\beta} r} g^{r \overline{\delta}} \chi_{\alpha \overline{\delta}}) + \sigma^2 \triangle_g v \right) \; {{\omega_{\varphi}^n}\over {n!}}  \\
& = &  \displaystyle \int_V \left( 2 \sum_{\alpha,\beta,r,\delta=1}^n (g^{\alpha \overline{\beta}} \sigma_{,\overline{\beta}} v_{,r} g^{r \overline{\delta}} \chi_{\alpha \overline{\delta}} + g^{\alpha \overline{\beta}} \sigma v_{,r} g^{r \overline{\delta}} \chi_{\alpha \overline{\delta},\overline{\beta}} )  - 2 \sigma \sigma_{,\overline{\delta}} v_{,r} g^{\overline{\delta} r} \right){{\omega_{\varphi}^n}\over {n!}} \\
& = &  \displaystyle \int_V \left( 2 \sum_{\alpha,\beta,r,\delta=1}^n (g^{\alpha \overline{\beta}} \sigma_{,\overline{\beta}} v_{,r} g^{r \overline{\delta}} \chi_{\alpha \overline{\delta}} + g^{\alpha \overline{\beta}} \sigma v_{,r} g^{r \overline{\delta}} \chi_{\alpha \overline{\beta},\overline{\delta}} )  - 2 \sigma \sigma_{,\overline{\delta}} v_{,r} g^{\overline{\delta} r} \right){{\omega_{\varphi}^n}\over {n!}}  \\
& = &  2\;\displaystyle \int_V \left(
\sum_{\alpha,\beta,r,\delta=1}^n g^{\alpha \overline{\beta}}
\sigma_{,\overline{\beta}} v_{,r} g^{r \overline{\delta}}
\chi_{\alpha \overline{\delta}} \right)
{{\omega_{\varphi}^n}\over {n!}} \\ &  = & - 2\;\displaystyle
\int_V \; \sum_{\alpha,\beta,r,\delta=1}^n \left( g^{\alpha
\overline{\beta}} \sigma_{,\overline{\beta}}\chi_{\alpha
\overline{\delta}}  g^{r \overline{\delta}} \right)_{,r}  \;
v\;{{\omega_{\varphi}^n}\over {n!}} .
\end{array}
\]
The Euler-Lagrange equation for functional $E$ is:
\begin{equation}
  \sum_{\alpha,\beta,r,\delta=1}^n \left( g^{\alpha \overline{\beta}} \sigma_{,\overline{\beta}}\chi_{\alpha \overline{\delta}}
    g^{r \overline{\delta}} \right)_{,r} = 0.
\label{eq:euler3}
\end{equation}
The left hand side is a divergence form and the equation holds on
the manifold without boundary. Multiple $\sigma$ in both sides of
the equation (\ref{eq:euler3}) and integrating over the entire
manifold, we have
\[\begin{array}{lcl} 0 & = & \displaystyle \int_V \;\left( g^{\alpha \overline{\beta}} \sigma_{,\overline{\beta}}\chi_{\alpha \overline{\delta}}
    g^{r \overline{\delta}} \right)_{,r} \sigma
    {{\omega_{\varphi}^n}\over {n!}} \\ & = &
  - \displaystyle \int_V \;  \sum_{\alpha,\beta,r,\delta=1}^n \left( g^{\alpha \overline{\beta}} \sigma_{,r}\chi_{\alpha \overline{\delta}}\;  g^{r \overline{\delta}} \right)  \sigma_{,\overline{\beta}} {{\omega_{\varphi}^n}\over {n!}} .
\end{array}
\]
Therefore, $\sigma = const$ in the manifold $V$ since $\chi$ is a
strictly positive (1,1) form.

 Now, along the $J-$flow, we have
${{\partial \varphi}\over {\partial t}} = c - tr_g \chi = c - \sigma.\;$ Thus,
\begin{eqnarray}
  {{d \; E}\over {d\,t}} & = &  2\;\displaystyle \int_V \left(  \sum_{\alpha,\beta,r,\delta=1}^n g^{\alpha \overline{\beta}}
  \sigma_{,\overline{\beta}} (c- \sigma)_{,r} g^{r \overline{\delta}} \chi_{\alpha \overline{\delta}} \right){{\omega_{\varphi}^n}\over {n!}} \nonumber \\
  & = &
 - 2\;\displaystyle \int_V \left(  \sum_{\alpha,\beta,r,\delta=1}^n g^{\alpha \overline{\beta}} \sigma_{,\overline{\beta}} \sigma_{,r} g^{r \overline{\delta}} \chi_{\alpha \overline{\delta}} \right){{\omega_{\varphi}^n}\over {n!}}
 \nonumber \\ & \leq & 0.
\label{eq:gradE}
\end{eqnarray}
The equality holds unless $\sigma $ is a constant. Thus $E$ is
strictly decreasing under the $J-$flow.
\end{proof}

\begin{prop} Any critical point of $E$ is a local minimizer.
\end{prop}
\begin{proof} Let $u,v$ be the tangential vertors in $V.\;$ Recall the first variation of $E$ as:
\[
\delta_v \; E(g) =  2\;\displaystyle \int_V \left(  \sum_{\alpha,\beta,r,\delta=1}^n g^{\alpha \overline{\beta}} \sigma_{,\overline{\beta}}
 v_{,r} g^{r \overline{\delta}} \chi_{\alpha \overline{\delta}} \right)\;{{\omega_{\varphi}^n}\over {n!}} .
\]
Next we want to compute the second variation of $E:$
\[
\begin{array} {lcl} & & \delta_u \; \delta_v \; E(g) \\
& = & - 2\; \displaystyle \int_V \left(
\sum_{\alpha,\beta,r,\delta=1}^n g^{\alpha \overline{\beta}}
\left( \displaystyle \sum_{i,j,k,l=1}^n\;g^{i \overline{j} }
u_{,\overline{j} k} g^{k \overline{l}} \chi_{k
\overline{l}}\right)_{,\overline{\beta}} v_{,r} g^{r
\overline{\delta}} \chi_{\alpha \overline{\delta}} \right)
{{\omega_{\varphi}^n}\over {n!}} \\ & = & 2\; \displaystyle
\int_V \left(  \sum_{\alpha,\beta,r,\delta,i,j,k,l=1}^n\;
g^{\alpha \overline{\beta}}\;g^{i \overline{j} } u_{,\overline{j}
k} g^{k \overline{l}} \chi_{k \overline{l}}\;
v_{,r\overline{\beta}} g^{r \overline{\delta}}
\chi_{\alpha \overline{\delta}} \right)  {{\omega_{\varphi}^n}\over {n!}} \\
 & & \qquad \qquad   + 2\; \displaystyle \int_V \left(  \sum_{\alpha,\beta,r,\delta,i,j,k,l=1}^n
 \;g^{i \overline{j} } u_{,\overline{j} k} g^{k \overline{l}} \chi_{k \overline{l}}\; v_{,r} g^{r \overline{\delta}}
  \;  g^{\alpha \overline{\beta}}\;\chi_{\alpha \overline{\delta},\overline{\beta}} \right){{\omega_{\varphi}^n}\over {n!}} \\
& = & 2\; \displaystyle \int_V \left(  \sum_{\alpha,\beta,r,\delta,i,j,k,l=1}^n\; g^{\alpha \overline{\beta}}\;g^{i \overline{j} }
 u_{,\overline{j} k} g^{k \overline{l}} \chi_{k \overline{l}}\; v_{,r\overline{\beta}} g^{r \overline{\delta}}
  \chi_{\alpha \overline{\delta}} \right) {{\omega_{\varphi}^n}\over {n!}}  \\
 & & \qquad \qquad \qquad   + 2\; \displaystyle \int_V \left(  \sum_{\alpha,\beta,r,\delta,i,j,k,l=1}^n\;g^{i \overline{j} } u_{,\overline{j} k} g^{k \overline{l}} \chi_{k \overline{l}}\; v_{,r} g^{r \overline{\delta}} \;  \sigma_{, \overline{\delta}} \right){{\omega_{\varphi}^n}\over {n!}} \\
& = & 2\; \displaystyle \int_V \left(
\sum_{\alpha,\beta,r,\delta,i,j,k,l=1}^n\; g^{\alpha
\overline{\beta}}\;g^{i \overline{j} } u_{,\overline{j} k} g^{k
\overline{l}} \chi_{k \overline{l}}\; v_{,r\overline{\beta}} g^{r
\overline{\delta}} \chi_{\alpha \overline{\delta}} \right)
{{\omega_{\varphi}^n}\over {n!}} > 0.

\end{array}
\]
The last equality holds since $\sigma = const$ at the critical
point. \end{proof}
\begin{prop} Along the $J-$ flow,  any smooth curve in $\cal H$ strictly decreases. Moreover, the distance between any two
points decreases as well.
\end{prop}
\begin{proof} Suppose $\varphi(s): [0,1] \rightarrow \cal H$ is a smooth
curve in the space of K\"ahler potentials. Now consider the
energy of this curve as
\[  En = \displaystyle \int_{s=0}^1\;\displaystyle \int_V \; \left({{\partial \varphi} \over {\partial s}}
\right)^2 \; {{\omega_{\varphi}^n}\over {n!}}\;d\,s . \] Under the
$J$ flow, suppose that the energy of the evolved curve at time
$t>0$ is  $En(t).\;$ We want to show that the energy is strictly
decreasing under this flow (All of the derivatives, norm and
integration are taken with respect to $g $ in the following
calculation):
\[\begin{array}{lcl}
  & &  {{d\, En(t) }\over {d\, t}}\\
   & = &  2 \;\displaystyle \int_{s=0}^1 \; \displaystyle\; \int_V \;
   \left( {{\partial \varphi} \over {\partial s}} {{\partial^2 \varphi} \over {\partial s \partial t }}
   + ({{\partial \varphi} \over {\partial s}})^2 \triangle_g {{\partial \varphi} \over {\partial t}} \right)
    {{\omega_{\varphi}^n}\over {n!}} \;d\,s\\
  & = & - 2 \;\displaystyle \int_{s=0}^1 \; \displaystyle\; \int_V \;\left( {{\partial \varphi} \over {\partial s}} {{\partial \sigma}
\over {\partial s}} + ({{\partial \varphi} \over {\partial s}})^2
\triangle_g \sigma \right) \;{{\omega_{\varphi}^n}\over
{n!}}\;d\,s \\ &  = &   2 \;\displaystyle \int_{s=0}^1 \;
\displaystyle\; \int_V \;{{\partial \varphi}\over {\partial
s}}\cdot  \left( \displaystyle \sum_{\alpha,\beta,r,\delta=1}^n
g^{\alpha \overline{\beta}} ({{\partial \varphi}\over {\partial
s}})_{,\overline{\beta} r} g^{r \overline{\delta}} \chi_{\alpha
\overline{\delta}} + ( {{\partial \varphi} \over {\partial
s}})_{,\alpha} \sigma_{,\overline{\beta}}
 g^{\alpha \overline{\beta}}\right) {{\omega_{\varphi}^n}\over {n!}} \;d\,s \\
& = & - 2 \;\displaystyle \int_{s=0}^1 \; \displaystyle\; \int_V
\; \displaystyle \sum_{\alpha,\beta,r,\delta=1}^n \left(g^{\alpha
\overline{\beta}} ( {{\partial \varphi}\over {\partial
s}})_{,\overline{\beta}}  \; g^{\alpha \overline{\beta}}
({{\partial \varphi}\over {\partial s}})_{, r} g^{r
\overline{\delta}} \chi_{\alpha \overline{\delta}} \right.\\
 & & \qquad + \left.{{\partial
\varphi}\over {\partial s}} \cdot   g^{\alpha \overline{\beta}}
({{\partial \varphi}\over {\partial s}})_{, r} g^{r
\overline{\delta}} \chi_{\alpha
\overline{\delta},\overline{\beta}}\right)
{{\omega_{\varphi}^n}\over {n!}} \;d\,s
+ 2 \;\displaystyle \int_{s=0}^1 \; \displaystyle\; \int_V \; {{\partial \varphi} \over {\partial s}}
 ( {{\partial \varphi} \over {\partial s}})_{,\alpha}  \sigma_{,\overline{\beta}} g^{\alpha \overline{\beta}} {{\omega_{\varphi}^n}\over {n!}} \;d\,s\\
& = & - 2 \;\displaystyle \int_{s=0}^1 \; \displaystyle\; \int_V \; \displaystyle \sum_{\alpha,\beta,r,\delta=1}^n \left(g^{\alpha \overline{\beta}}
 ( {{\partial \varphi}\over {\partial s}})_{,\overline{\beta}}  \; g^{\alpha \overline{\beta}} ({{\partial \varphi}\over {\partial s}})_{, r} g^{r \overline{\delta}} \chi_{\alpha
 \overline{\delta}}\right. \\ & & \qquad
 \left. + {{\partial \varphi}\over {\partial s}} \cdot   ({{\partial \varphi}\over {\partial s}})_{, r} g^{r \overline{\delta}} \sigma_{,\overline{\delta}}\right){{\omega_{\varphi}^n}\over {n!}} \;d\,s
+ 2 \;\displaystyle \int_{s=0}^1 \; \displaystyle\; \int_V \;
{{\partial \varphi} \over {\partial s}} ( {{\partial \varphi}
\over {\partial s}})_{,\alpha}  \sigma_{,\overline{\beta}}
g^{\alpha \overline{\beta}} \;{{\omega_{\varphi}^n}\over {n!}}
\;d\,s\\ & = & - 2 \;\displaystyle \int_{s=0}^1 \;
\displaystyle\; \int_V \; \displaystyle
\sum_{\alpha,\beta,r,\delta=1}^n \left(g^{\alpha
\overline{\beta}} ( {{\partial \varphi}\over {\partial
s}})_{,\overline{\beta}}  \; g^{\alpha \overline{\beta}}
({{\partial \varphi}\over {\partial s}})_{, r} g^{r
\overline{\delta}} \chi_{\alpha \overline{\delta}}\right)
{{\omega_{\varphi}^n}\over {n!}} \;d\,s \leq 0.
\end{array}
\]
The equality holds unless ${{\partial \varphi}\over {\partial s}} $ is a constant in $ V\times [0,1]\;$ or the geodesic is trivial. \\

If the $J$ flow exists for all the time for any initially smooth
metric, then this will imply that $J$ flow decreases distance
between any two points in $\cal H.\;$ We will prove the long term
existence of the $J$ flow in the next section.
\end{proof}

There are surprising similarities between these two group of
functionals: the first group is the Mabuchi energy and the Calabi
energy, both are well known, but perhaps not well understood. The
second group is our $J$ functional and its norm $E.\;$ Here $J$
plays the "role" of the Mabuchi energy and $E$ plays the role of
the Calabi energy.  Calabi first proved in \cite{calabi82} that
the critical point of the Calabi energy is a local minimizer. We
shall compare this to Proposition 2.4 aforementioned. In
\cite{chen992}, we showed that the Calabi flow (gradient flow of
The Mabuchi energy ) decreases the length of any smooth curves in
the space of K\"ahler metrics. We shall compare this to
Proposition 2.5 aforementioned. In \cite{chen991}, we prove that
the critical point of the Calabi energy is unique if the first
Chern class is negative; and in general it is conjectured that
critical point of the Calabi energy is unique in each K\"ahler
class. We shall compare this to Proposition 2.1 aforementioned.
The list of similarities can go on and on.  The critical point of
the Calabi energy is well known. The critical point of $E$ is not
so well-known and is also not clear about its significance in
geometry. But the similarity between $E$ and the Calabi energy
makes one thing  clear: to study the critical point of $E$ or its
gradient flow, is amount to study a junior version of extremal
metrics or the Calabi flow. The insight we will learn and the
technique we will develop from the studying of the Euler-Lagrange
equation for the functional $J$ must be helpful to the studying of
the extremal metric and the Calabi flow.

\section{$C^2$ estimate of heat flow depending time $t$}
\begin{lem} $\sigma > 0$ is bounded from above and below along
the gradient flow of $J.\;$
\end{lem}
\begin{proof} Taking second derivatives with respect time, we have
\[
 {{\partial^2 \varphi}\over {\partial t^2}} =
 \displaystyle \sum_{\alpha,\beta,r,\delta=1}^n g^{\alpha \overline{\beta}} \left({{\partial \varphi}\over {\partial t}}\right)_{,\overline{\beta} r}
 g^{r \overline{\delta}} \chi_{\alpha \overline{\delta}}.
\]
By the ordinary maximum principle for parabolic equation, we have
$\displaystyle \max_{V} {{\partial \varphi}\over {\partial t}}$
is decreasing while $\displaystyle \min_{V} {{\partial
\varphi}\over {\partial t}} $ is increasing as $t $ increases.
Since ${{\partial \varphi}\over {\partial t}} = c - \sigma,$ we
then prove this lemma. \end{proof}

\begin{cor} Under  the  $J-$flow  the evolved metric is
strictly bounded from below.
\end{cor}

\begin{theo} There exists a const $C(t)$ depends on $t$
only, such that  $\omega_0 + \partial \overline{\partial} \varphi
\leq C \cdot \omega_0.\;$
\end{theo}
\begin{proof} Consider the heat flow equation:
\[ {{\partial \varphi}\over {\partial t}} =
c - \displaystyle \sum_{\alpha, \beta=1}^n\; g^{\alpha\;
\overline{\beta}} \;\chi_{\alpha \overline{\beta}}.\]
 In the
following calculation, all indices are running from $1$ to $n.\;$
We adopt Einstein's convention that the repeated indices  mean
summation from $1$ to $n.\;$Taking first derivatives on both side
of this equation.
\[
  {{\partial \varphi_{,i}}\over {\partial t}} = g^{\alpha \overline{r}} g_{\overline{r} \delta,i} g^{\delta \overline{\beta}} \chi_{\alpha \overline{\beta}} - g^{\alpha \overline{\beta}} \chi_{\alpha \overline{\beta},i}.
\]
Next we take second derivatives on both sides of the above
equation. We have
\begin{equation}
\begin{array}{lcl}  {{\partial \varphi_{i \overline{j}}} \over {\partial t}}
& = & g^{\alpha \overline{r}} g_{\overline{r} \delta,i
\overline{j}} g^{\delta \overline{\beta}} \chi_{\alpha
\overline{\beta}} - g^{\alpha \overline{a}} g_{\overline{a}
b,\overline{j}} g^{b \overline{r}} g_{\overline{r} \delta,i}
g^{\delta \overline{\beta}} \chi_{\alpha \overline{\beta}} \nonumber \\
& & \qquad - g^{\alpha \overline{r}}
g_{\overline{r} \delta,i} g^{\delta \overline{a}} g_{\overline{a} b,\overline{j}} g^{b \overline{\beta}} \chi_{\alpha \overline{\beta}}
  - {g^{\alpha \overline{\beta}}}_{\overline{j}} \chi_{\alpha
\overline{\beta}, i}  - g^{\alpha \overline{\beta}} \chi_{\alpha
\overline{\beta}, i \overline{j}}. \label{eq:C2estimate1}
\end{array}
\end{equation}
Indices after comma represent partial derivatives in above two
equations. Next we multiple the  equation (\ref{eq:C2estimate1})
with $\chi^{i \overline{j}}$ and summarize everything. Set
\[ F =
\chi^{i \overline{j}} g_{i\overline{j}} =  \chi^{i \overline{j}}
({g_0}_{i\overline{j}} + \varphi_{,i \overline{j}}).\]

 Consider
all differentials as covariant differentials with respect to
metric $\chi.\;$ Denote the bisectional curvature of $\chi$ as
$R(\chi).\;$We then have
\[
\begin{array}{lcl} {{\partial F}\over {\partial t}} & = &  g^{\alpha \overline{r}} g_{\overline{r} \delta,i \overline{j}}
g^{\delta \overline{\beta}} \chi_{\alpha \overline{\beta}} \chi^{i \overline{j}} -2   g^{\alpha \overline{a}}
 g_{\overline{a} b,\overline{j}} g^{b \overline{r}} g_{\overline{r} \delta,i} g^{\delta \overline{\beta}}
 \chi_{\alpha \overline{\beta}} \chi^{i \overline {j}}\\
 & = & g^{\alpha \overline{r}} g_{\overline{r} i,\delta \overline{j}} g^{\delta \overline{\beta}} \chi_{\alpha \overline{\beta}}
  \chi^{i \overline{j}} -2   g^{\alpha \overline{a}} g_{\overline{a} b,\overline{j}} g^{b \overline{r}} g_{\overline{r} \delta,i}
  g^{\delta \overline{\beta}} \chi_{\alpha \overline{\beta}} \chi^{i \overline {j}}\\
& = & g^{\alpha \overline{r}} \left(g_{\overline{r} i,\overline{j} \delta}
 + R(\chi)_{i \overline{p} \delta \overline{j}} g_{p \overline{r}}
 - R(\chi)_{p \overline{r} \delta \overline{j}} g_{i \overline{p}} \right)
 g^{\delta \overline{\beta}} \chi_{\alpha \overline{\beta}} \chi^{i \overline{j}}
\\
& & \qquad \qquad \qquad \qquad -2   g^{\alpha \overline{a}}
g_{\overline{a} b,\overline{j}} g^{b \overline{r}}
g_{\overline{r} \delta,i}
 g^{\delta \overline{\beta}} \chi_{\alpha \overline{\beta}} \chi^{i \overline {j}}\\
 & = &   g^{\alpha \overline{r}}g_{ i\overline{j}, \overline{r}\delta}g^{\delta \overline{\beta}} \chi_{\alpha \overline{\beta}}
 \chi^{i \overline{j}} +
g^{\alpha \overline{r}} \left( R(\chi)_{i \overline{p} \delta \overline{j}} g_{p \overline{r}}
- R(\chi)_{p \overline{r} \delta \overline{j}} g_{i \overline{p}} \right) g^{\delta \overline{\beta}} \chi_{\alpha \overline{\beta}}
 \chi^{i \overline{j}}\\
 & & \qquad \qquad \qquad \qquad
 -2   g^{\alpha \overline{a}} g_{\overline{a} b,\overline{j}} g^{b \overline{r}} g_{\overline{r} \delta,i} g^{\delta \overline{\beta}}
  \chi_{\alpha \overline{\beta}} \chi^{i \overline {j}}. \end{array}
\]
Since $\chi$ is a fixed K\"ahler metric, then there exists a
constant $C$ such that
\[\mid R(x)_{i \overline{j} k \overline{l}
}\mid \leq C \cdot (\chi_{i \overline{j} } \chi_{k \overline{l}} +
\chi_{i \overline{l} } \chi_{k \overline{j}}).
\]
Also, we can
choose a coordinate such that $\chi_{i \overline{j}} = \delta_{i
\overline{j}}.\; $ Thus,
\[\begin{array}{lcl} & & g^{\alpha \overline{r}} \left( R(\chi)_{i \overline{p} \delta \overline{j}} g_{p \overline{r}}
- R(\chi)_{p \overline{r} \delta \overline{j}} g_{i \overline{p}}
\right) g^{\delta \overline{\beta}} \chi_{\alpha
\overline{\beta}} \chi^{i \overline{j}} \\ &  & \leq  C \cdot
g^{\alpha \overline{r}}  \left( ( \delta_{i \overline{p}}
\delta_{\delta \overline{j}} +  \delta_{i \overline{j}}
\delta_{\delta \overline{p}}) g_{p \overline{r}} +   ( \delta_{p
\overline{r} }\delta_{\delta \overline{j}} +  \delta_{p
\overline{j}} \delta_{\delta \overline{r}} ) g_{i \overline{p}}
\right)g^{\delta \overline{\beta}} \delta_{\alpha
\overline{\beta}} \delta^{i \overline{j}} \\&  & \leq  C \cdot
g^{\alpha \overline{\alpha}} = C \;\sigma.
\end{array}
\]
Thus, we have
\[
\begin{array}{lcl}
{{\partial F} \over {\partial t}}
& = & \tilde{\triangle} F - 2   g^{\alpha \overline{a}} g_{\overline{a} b,\overline{j}} g^{b \overline{r}} g_{\overline{r} \delta,i} g^{\delta \overline{\beta}} \chi_{\alpha \overline{\beta}} \chi^{i \overline {j}} + C\cdot \sigma.
\end{array}
\]
Where $\tilde{\triangle} f = g^{\alpha \overline{r}}
f_{,\overline{r}\delta} g^{\delta \overline{\beta}} \chi_{\alpha
\overline{\beta}} $ for any smooth function $f.\;$  From here, we
quickly deduce that $F$ is bounded from above since $\sigma $ is
uniformly bounded from above.
\end{proof}

Recall the estimate of equation ~(\ref{eq:gradE}), we have
\[
\begin{array}{lcl}
  \displaystyle \int_0^{\infty}\; {{d \; E}\over {d\,t}}\; d\;t & = &   - 2 \;\displaystyle \int_0^{\infty} \;\displaystyle
  \int_V \left(  \sum_{\alpha,\beta,r,\delta=1}^n g^{\alpha \overline{\beta}} \sigma_{,\overline{\beta}} \sigma_{,r}
   g^{r \overline{\delta}} \chi_{\alpha \overline{\delta}} \right) \;{{\omega_{\varphi}^n}\over {n!}}
   \;d\;t\\
 & = & E(\infty) - E(0) \geq C.
 \end{array}
\]
In other words,
\[
\;\displaystyle \int_0^{\infty} \;\displaystyle \int_V \left(
\sum_{\alpha,\beta,r,\delta=1}^n g^{\alpha \overline{\beta}}
\sigma_{,\overline{\beta}} \sigma_{,r} g^{r \overline{\delta}}
\chi_{\alpha \overline{\delta}} \right)
\;{{\omega_{\varphi}^n}\over {n!}}  \;d\;t \leq C.
\]
There exists a subsequence of $t_i \rightarrow \infty $
such that
\[
\displaystyle \int_V \left(  \sum_{\alpha,\beta,r,\delta=1}^n
g^{\alpha \overline{\beta}} \sigma_{,\overline{\beta}}
\sigma_{,r} g^{r \overline{\delta}} \chi_{\alpha
\overline{\delta}} \right) \;{{\omega_{\varphi}^n}\over {n!}}
\mid_{t=t_i} \rightarrow 0.
\]
This last expression, shall suggest that $\sigma\mid_{t=t_i}
\rightarrow c$
 in some  sense for some constant $c$.\\

Next we return to prove the first part of Theorem 0.1.
\begin{proof} Following
from the interior estimate by Evans and Krylov, we can imply
$C^{2,\alpha}$ estimate for any finite time $t.\;$ Then standard
theory of the parabolic equation  will  imply that $g$ is
$C^{\infty}$ at any finite time. Thus the flow exists for long
time.  Then proposition 2.5
 implies that $J$ flow decreases
distance between any two points in $\cal H.\;$
\end{proof}

\section{Uniform $C^2$ estimate for heat flow for manifolds
with semi-positive definite curvature tensors.}
In this section, we assume that the bisectional curvature of $\chi$ is
non-negative, we want to show that there exists a uniform bound on the
second derivatives of $\varphi.\;$
\begin{theo} If the bisectional curvature of $\chi$ is
non-negative, then there exists a uniform bound on the
second derivatives of $\varphi.\;$
\end{theo}
\begin{proof} Following the calculation in Section 3, we
have (equation (\ref{eq:C2estimate1})) :
\[
\begin{array} {lcl}  {{\partial \varphi_{,i \overline{j}}} \over {\partial t}}
& = & g^{\alpha \overline{r}} g_{\overline{r} \delta,i
\overline{j}} g^{\delta \overline{\beta}} \chi_{\alpha
\overline{\beta}} - g^{\alpha \overline{a}} g_{\overline{a}
b,\overline{j}} g^{b \overline{r}} g_{\overline{r} \delta,i}
g^{\delta \overline{\beta}} \chi_{\alpha \overline{\beta}} \\ & &
\qquad - g^{\alpha \overline{r}} g_{\overline{r} \delta,i}
g^{\delta \overline{a}} g_{\overline{a} b,\overline{j}} g^{b
\overline{\beta}} \chi_{\alpha \overline{\beta}}  - {g^{\alpha
\overline{\beta}}}_{\overline{j}} \chi_{\alpha \overline{\beta},
i}  - g^{\alpha \overline{\beta}} \chi_{\alpha \overline{\beta},
i \overline{j}}.
\end{array}
\]
Simplifying this equation by using  covariant derivatives in
terms of the metric $\chi,$  we have
\begin{eqnarray}
  {{\partial \varphi_{,i \overline{j}}} \over {\partial t}} & =  &  g^{\alpha \overline{r}} g_{\overline{r} \delta,i \overline{j}} g^{\delta \overline{\beta}} \chi_{\alpha \overline{\beta}}
   -  g^{\alpha \overline{a}} g_{\overline{a} b,\overline{j}} g^{b \overline{r}} g_{\overline{r} \delta,i} g^{\delta \overline{\beta}} \chi_{\alpha \overline{\beta}}
  \nonumber \\ & & \qquad \qquad \qquad \qquad- g^{\alpha \overline{r}}
g_{\overline{r} \delta,i} g^{\delta \overline{a}} g_{\overline{a} b,\overline{j}} g^{b \overline{\beta}} \chi_{\alpha \overline{\beta}}
\label{eq:c2estimate2}
\end{eqnarray}
Define an auxiliary tensor $T_{i \overline{j}}$ as

\[T_{i \overline{j} } = g_{i \overline{j}} - C_0 \cdot \chi_{i \overline{j}}
=  {g_0}_{i \overline{j}}  + \varphi_{i \overline{j}} - C_0 \cdot
\chi_{i \overline{j}}.\] Choose $C_0$ big enough so that
$T_{i\overline{j}} < 0$ as a tensor at time $t=0.\;$ \\

\noindent {\bf Claim}: $T_{i \bar j}<0$ is preserved under this
gradient flow. \\

From equation (\ref{eq:c2estimate2}), we have

\begin{equation}
\begin{array}{lcl}  {{\partial T_{,i \overline{j}}} \over {\partial t}}
& = &  g^{\alpha \overline{r}} T_{\overline{r} \delta,i \overline{j}} g^{\delta \overline{\beta}} \chi_{\alpha \overline{\beta}} -  g^{\alpha \overline{a}} g_{\overline{a} b,\overline{j}} g^{b \overline{r}} g_{\overline{r} \delta,i} g^{\delta \overline{\beta}} \chi_{\alpha \overline{\beta}} - g^{\alpha \overline{r}}
g_{\overline{r} \delta,i} g^{\delta \overline{a}} g_{\overline{a} b,\overline{j}} g^{b \overline{\beta}} \chi_{\alpha \overline{\beta}} \nonumber \\
& = &  g^{\alpha \overline{r}} T_{\overline{r} i, \delta \overline{j}} g^{\delta \overline{\beta}} \chi_{\alpha \overline{\beta}} -  g^{\alpha \overline{a}} g_{\overline{a} b,\overline{j}} g^{b \overline{r}} g_{\overline{r} \delta,i} g^{\delta \overline{\beta}} \chi_{\alpha \overline{\beta}} - g^{\alpha \overline{r}}
g_{\overline{r} \delta,i} g^{\delta \overline{a}} g_{\overline{a} b,\overline{j}} g^{b \overline{\beta}} \chi_{\alpha \overline{\beta}} \nonumber
\\ & = &  g^{\alpha \overline{r}} \left(T_{\overline{r}
i,\overline{j} \delta} + R(\chi)_{i \overline{p} \delta
\overline{j}} T_{p \overline{r}} - R(\chi)_{p \overline{r} \delta
\overline{j}} T_{i \overline{p}} \right) g^{\delta
\overline{\beta}} \chi_{\alpha \overline{\beta}}\\
& & \qquad \qquad -  g^{\alpha \overline{a}} g_{\overline{a}
b,\overline{j}} g^{b \overline{r}} g_{\overline{r} \delta,i}
g^{\delta \overline{\beta}} \chi_{\alpha \overline{\beta}} -
g^{\alpha \overline{r}}
g_{\overline{r} \delta,i} g^{\delta \overline{a}} g_{\overline{a} b,\overline{j}} g^{b \overline{\beta}} \chi_{\alpha \overline{\beta}} \nonumber \\
& = &  g^{\alpha \overline{r}} \left(T_{
i\overline{j},\overline{r} \delta} + R(\chi)_{i \overline{p}
\delta \overline{j}} T_{p \overline{r}} - R(\chi)_{p \overline{r}
\delta \overline{j}} T_{i \overline{p}} \right) g^{\delta
\overline{\beta}} \chi_{\alpha \overline{\beta}} \\ & &   \qquad
\qquad - g^{\alpha \overline{a}} g_{\overline{a} b,\overline{j}}
g^{b \overline{r}} g_{\overline{r} \delta,i} g^{\delta
\overline{\beta}} \chi_{\alpha \overline{\beta}} - g^{\alpha
\overline{r}} g_{\overline{r} \delta,i} g^{\delta \overline{a}}
g_{\overline{a} b,\overline{j}} g^{b \overline{\beta}}
\chi_{\alpha \overline{\beta}}.
 \label{eq:c2estimates3}
\end{array}
\end{equation}
Next we want to apply Hamilton's maximal principal for tensors.
Since $T_{i\overline{j}} < 0$ at $t=0,$ that there is a first
time $t=t_0 > 0 $ and a point $O,$ where $T$ has a degenerate
direction in $T_O V.\;$  We assume that this direction is $\xi
(O) = (\xi^{1}, \xi^2,\cdots \xi^n).\; $ Parallel transport this
vector along a small neighborhood $\cal O $ of $O$ by metric
$\chi.\;$ By definition, we have
\[
  T_{i \overline{j}} \xi^i  \overline{\xi}^j < 0, \;\; \forall\; t <
  t_0,
\]
and at $t=t_0$ we have
\[
  T_{i \overline{j}} \xi^i  \overline{\xi}^j(O) = 0,\qquad {\rm
  and}\qquad
\;T_{i \overline{j}} \leq 0 \qquad {\rm in} \; \cal O.
\]
In particular, at $t=t_0$ and at point $O$, we have
$T_{i\overline{j}} \xi^{i} = T_{i \overline{j}}
\overline{\xi^{j}} = 0.\;$ Now plugging everything into the
equation (\ref{eq:c2estimates3}), we have
\[\begin{array}{lcl}
& &  {{\partial (T_{i \overline{j}}} \xi^i \overline{\xi}^j) \over
{\partial t}} \\
& = &  g^{\alpha \overline{r}} T_{ i\overline{j},\overline{r}
\delta}\xi^i \overline{\xi}^j g^{\delta \overline{\beta}}
\chi_{\alpha \overline{\beta}}  +
  g^{\alpha \overline{r}} \left( R(\chi)_{i \overline{p} \delta \overline{j}} T_{p \overline{r}} - R(\chi)_{p \overline{r} \delta \overline{j}} T_{i \overline{p}} \right)
   g^{\delta \overline{\beta}} \chi_{\alpha \overline{\beta}}  \xi^i  \overline{\xi}^j \\ &   & \qquad   \qquad  \qquad  \qquad -  g^{\alpha \overline{a}} g_{\overline{a} b,\overline{j}} g^{b \overline{r}} g_{\overline{r} \delta,i} g^{\delta \overline{\beta}} \chi_{\alpha \overline{\beta}}  \xi^i  \overline{\xi}^j - g^{\alpha \overline{r}}
g_{\overline{r} \delta,i} g^{\delta \overline{a}} g_{\overline{a} b,\overline{j}} g^{b \overline{\beta}} \chi_{\alpha \overline{\beta}}  \xi^i  \overline{\xi}^j \\
& = & \tilde{\triangle} (T_{i \overline{j}} \xi^i
\overline{\xi}^j) +
 g^{\alpha \overline{r}}  R(\chi)_{i \overline{p} \delta \overline{j}} \xi^i \overline{\xi}^j T_{p \overline{r}}  g^{\delta \overline{\beta}} \chi_{\alpha \overline{\beta}} - g^{\alpha \overline{r}} R(\chi)_{p \overline{r} \delta \overline{j}} T_{i \overline{p}} \xi^i g^{\delta \overline{\beta}} \chi_{\alpha \overline{\beta}}    \overline{\xi}^j \\
&  & \qquad   \qquad  \qquad  \qquad -  g^{\alpha \overline{a}}
g_{\overline{a} b,\overline{j}} g^{b \overline{r}}
g_{\overline{r} \delta,i} g^{\delta \overline{\beta}}
\chi_{\alpha \overline{\beta}}  \xi^i  \overline{\xi}^j -
g^{\alpha \overline{r}} g_{\overline{r} \delta,i} g^{\delta
\overline{a}} g_{\overline{a} b,\overline{j}} g^{b
\overline{\beta}} \chi_{\alpha \overline{\beta}}  \xi^i
\overline{\xi}^j.
\end{array}
\]
Now at the point $O$ and at time $t=t_0,$ by the standard maximum
principle, we have  $\tilde{\triangle} (T_{i \overline{j}} \xi^i
\overline{\xi^j}) \leq 0.\;$
 Moreover,
 \[ g^{\alpha \overline{r}} R(\chi)_{p \overline{r} \delta \overline{j}} T_{i \overline{p}} \xi^i g^{\delta \overline{\beta}}
  \chi_{\alpha \overline{\beta}}    \overline{\xi}^j (O)= 0
  \] and
  \[  g^{\alpha \overline{r}}  R(\chi)_{i \overline{p} \delta \overline{j}} \xi^i \overline{\xi}^j T_{p \overline{r}}
    g^{\delta \overline{\beta}} \chi_{\alpha \overline{\beta}} \leq
    0.
  \]
The last inequality holds since $R(\chi) $ is a non-negative
tensor while $T$ is a non-positive tensor. Thus
\[
  {{\partial (T_{i \overline{j}}} \xi^i \overline{\xi}^j) \over {\partial t}}(O) \leq 0.
\]
This implies that $T$ will remain non-positive. In other words,
\[
  g_{i\overline{j}} \leq C_0 \cdot \chi_{i \overline{j}}
\]
holds for all $t$ where the flow exists. Thus, all of the second
derivatives of $\varphi$ is bounded from above.
\end{proof}

Finally, we want to prove the second part of Theorem 0.1.
 \begin{proof} Again,
 following from the interior estimate by Evans and Krylov, we can
obtain a uniform $C^{2,\alpha}$ estimate for any finite time $t$
from the Theorem 4.1  above. Then standard elliptic regularity
theorem would imply that $g$ is $C^{\infty}$ at any finite time.
Thus the flow exists for long time. Since the $C^{2,\alpha}$
estimate is uniform (independent of time $t$), thus the flow
converges to a critical point of $J,\;$ at least by sequence. The
uniqueness of sequential limit is provided by the fact that $J$
is strictly convex.
\end{proof}


\end{document}